\documentclass[numreferences]{kluwer}
\usepackage{graphicx}
\font\ggo=eufb10
\font\bBB=msbm10
\font\bbB=msbm7
\def\bBC{\mbox{\bBB C}}
\def\bBc{\mbox{\bbB C}}

\def\bBr{\mbox{\bbB R}}

\def\ggG{\mbox{\ggo G}}
\def\ggB{\mbox{\ggo B}}

\def\ggD{\mbox{\ggo D}}
\begin{document}
\begin{article}
\begin{opening}
\title{Irreducible Representations of a Class of Current Algebras of Etingof and Frenkel\thanks{Submitted to Letters in Mathematical Physics}}
\author{O. \surname{Stoytchev}\email{ostoytchev@aubg.bg}}
\institute{American University in Bulgaria, 2700 Blagoevgrad, Bulgaria\\
and Institute For Nuclear Research, 1784 Sofia, Bulgaria}
\begin{abstract}
A class of  representations is described for the central extensions, found by Etingof and Frenkel \cite{Etingof}, of current algebras over Riemann surfaces. Their irreducibility is proved. The possibility/impossibility to obtain integrable representations within that class is discussed briefly.
\end{abstract}
\keywords{Representation Theory, Current Algebras, Central Extensions, Infinite Dimensional Lie Algebras and Groups, Riemann Surfaces}
\end{opening}
\begin{section}{Introduction}
 In this paper $\ggG$ will be a semisimple Lie algebra and $X$ a 
Riemannian manifold. The space $C_0^\infty(X,\ggG)\equiv\ggG^X$ of compactly
supported smooth functions with values in $\ggG$ has a natural structure of
(infinite dimensional) Lie algebra with respect to the pointwise Lie bracket
$$
[a,b](x):=[a(x),b(x)]\,,\qquad a,b\in\ggG^X\,.
$$
Let $G$ be the simply connected Lie group whose Lie algebra is $\ggG$.
By $G^X$ we shall denote the space of smooth maps from $X$ into $G$, which are
identity outside a compact subset of $X$. It becomes a topological
group with the obvious pointwise group operations induced from those on $G$ and
the natural $C^{\infty}$ topology. The pointwise exponential mapping
is a local homeomorphism from $\ggG^X$ into $G^X$.
We shall use the term current group, resp. current algebra for  $G^X$, resp. $\ggG^X$, or any central extension of these. \par
The representation theory of loop groups and loop algebras and their central extensions, i.e. the one-dimensional 
case $X=S^1$, has been developed extensively
during the past two decades. Numerous important results are contained in the well-known monographs
by Kac \cite{Kac} (at the algebra level) and Pressley and Segal \cite{Pressley}.  By contrast, the representation theory 
of current groups and current algebras over higher-dimensional
manifolds contains relatively few results. The two-dimensional case promises to reveal a rich mathematical 
structure, as noted in the paper by Etingof and Frenkel \cite{Etingof}. In that work the authors constructed and studied interesting central extensions of $\ggG^X$ (we call them for brevity Etingof--Frenkel current algebras) in the case when $X$ is a compact Riemann surface and $\ggG$ is a complex simple Lie algebra. These extensions are distinguished by the property that they can be integrated to  extensions by  tori of the corresponding current groups. \par
   In a short announcement \cite{Stoytchev} we outlined a construction of a class of irreducible representations of the Etingof-Frenkel current algebras. The only known (interesting) irreducible representations of current groups in dimension higher than one are variations of the so called energy representations, constructed by Albeverio and H{\o}egh-Krohn \cite{Albeverio}, Ismagilov \cite{Ismagilov} and Gelfand, Graev and Ver\v{s}ik \cite{{Gelfand},{Gelfand1}}. These are representations of the nonextended groups $G^X$ and the corresponding algebra representations are obtained from them. A modification, which we called the restricted energy representation, was found in  \cite{Donkov}. It is always irreducible in dimension two and higher, gives rise to an irreducible representation of the corresponding algebra and the latter allows a natural extension to a representation of certain central extensions of that algebra. Our initial scheme did not include directly representations of the Etingof--Frenkel current algebras. The representations of the latter, found in \cite{Stoytchev},  have an interesting property --- the nontrivial topology of $X$ manifest itself. This is not the case for any of the previously known representations. Namely, for those representations, if we take e.g. $X_1$ to be a 2-torus and $X_2$ to be a square, the representation space for $\ggG^{X_1}$ and $\ggG^{X_2}$ is the same and the algebra $\ggG^{X_1}$ acts on it as the
subalgebra of $\ggG^{X_2}$ of functions with periodic boundary conditions. By contrast, to construct irreducible representations of the Etingof--Frenkel current algebra for the torus, we need to tensor our ``old'' representation space with two irreducible highest (or lowest) weight modules over Kac--Moody algebras, which are restrictions of $\ggG^{X_1}$ to the two cycles, generating the first homology group of the torus.\par
  The next section is an exposition of results on representations of current algebras which are not of the Etingof--Frenkel type. The proof of Proposition 2 has not been published before. In Section 3 we explain how the Etingof--Frenkel current algebras are defined and construct irreducible representations of those. Our main result is the irreducibility result contained in Proposition 3. As far as we know, the representations that we have found are the only known irreducible representations of Etingof--Frenkel algebras.
\par
In the last section we give two examples of representations over the torus and over a hyperelliptic surface of genus 2. These examples as well as the general theory of holomorphic differentials over Riemann surfaces indicate that, unfortunately, one cannot hope to obtain integrable representations of the Etingof--Frenkel algebras using just integrable highest (lowest) weight modules over Kac--Moody algebras.
\end{section}
\begin{section}{Irreducible representations of $\ggG^X$ and some of its central extensions}
   In this section we state and prove some important preliminary results. The manifold $X$ will be an arbitrary compact two-dimensional Riemannian manifold and $\ggG$ will be a compact simple Lie algebra. (All results are actually valid for noncompact manifolds if we require our algebras to consist of functions with compact support.) 
 Consider the Hilbert space $L^2(X,dx;\ggG_{\bBc})$, where $d\mu (x)$ is the
Riemann-Lebesgue measure on the manifold $X$ and the scalar product is
\begin{equation}
(a,b):=\int_X\,\langle \bar{a}(x),b(x)\rangle\,d\mu (x) \label{scalar}
\end{equation}
with $\langle\ ,\ \rangle$ an invariant inner product on $\ggG$. The
representations which we consider are realized in the bosonic Fock space over
the one-particle Hilbert space $L^2(X,dx;\ggG_{\bBc})$. We denote by $B^*(a)$ and
$B(a)$ the standard creation and annihilation operators (of a state $a$) in the
Fock space and by $L_0(a)$ --- the second-quantized operator, corresponding
to the one-particle operator $b(x)\rightarrow [a(x),b(x)]$. On an $n$-particle
state these act as follows:
\begin{eqnarray}
\nonumber B^*(a)\,(B^*(a_1)B^*(a_2)...B^*(a_n)|0\rangle)& \!\!\!=\!\!\! &
B^*(a)B^*(a_1)B^*(a_2)...B^*(a_n)|0\rangle\,, \\                         
\nonumber B(a)\,(B^*(a_1)B^*(a_2)...B^*(a_n)|0\rangle)& \!\!\!=\!\!\! &
\sum_{i=1}^n\,(a,a_i)\,B^*(a_1)...B^*(a_{i-1})\\ 
\nonumber &\!\!\!\times\!\!\! &B^*(a_{i+1})...B^*(a_n)|0\rangle , \\ 
\nonumber L_0(a)\,(B^*(a_1)B^*(a_2)...B^*(a_n)|0\rangle)& \!\!\!=\!\!\! &
\sum_{i=1}^nB^*(a_1)...B^*(a_{i-1})B^*([a,a_i]) \\ 
\nonumber &\!\!\!\times\!\!\! & B^*(a_{i+1})...B^*(a_n)|0\rangle \,.
\end{eqnarray}
   Take a piecewise smooth real vector field $Z$ on $X$, which is nonzero
a.e.
(By piecewise smooth we mean that the coefficient functions of $Z$ in any
coordinate system are smooth and bounded together with all their derivatives
outside a closed subset of $X$ of zero measure, i.e.~we allow jumps.
In this way the
nontrivial topology of $X$ doesn't play any role at this point.) The vector
field $Z$ is the data which will determine the representation. Recall that
vector fields are operators on (smooth) functions on $X$.

All our irreducibility results will rely on the following proposition, proven in \cite{Donkov}:
\newtheorem{theoremdemo}{Proposition}
\begin{theoremdemo}
The operators
\begin{equation}
L^Z(a):=L_0(a)+B^*(Za)-B(Z{a})\,, \qquad a\in\ggG^X\,, \label{realrep}
\end{equation}
define a topologically irreducible unitary $($i.e.~antihermitian$)$ representation of the
algebra $\ggG^X$. For a fixed Riemannian metric, different vector fields determine
nonequivalent representations. The metric itself is unimportant since what
actually matters is the Riemann-Lebesgue measure, i.e.~the volume form, but
any local rescaling of the measure turns out to be equivalent $($in the sense
of representation theory$)$ to a rescaling of $Z$.
\end{theoremdemo}
We would like to point out that although the action of $L^Z(a)$ (Eq. (\ref{realrep})) is fairly simple, the proof of irreducibility (that we know) is difficult. It is analytical and is done at the group level. Note that in the case Dim~$X=1$ the representation (Eq. (\ref{realrep})) is known to be reducible \cite{Albeverio2}.\par
   It was noticed \cite{Donkov} that if one replaces the real vector field $Z$ in Eq. (\ref{realrep}) by a complex vector field $Y=Z_1+iZ_2$, then the operators $L^Y$ defined as in Eq. (\ref{realrep}) with $Y$ instead of $X$ give a representation of a central extension of $\ggG^X$. Namely, we have:
\begin{equation}
L^Y(a)L^Y(b)-L^Y(b)L^Y(a)=L^Y([a,b])+c(a,b)\,, \label{extended}
\end{equation} 
where
\begin{equation}
c(a,b)=-2i\left[(Z_1{b},Z_2a)-(Z_1{a},Z_2b)\right] \label{c}
\end{equation}
is  2-cocycle on $\ggG^X$, which is nontrivial unless it is zero \cite{Donkov}.
\begin{theoremdemo}
The operators $L^Y(a)$ defined as in {\rm Eq. (\ref{realrep})} with complex vector field $Y$ give a unitary (i.e. the operators are antihermitian) topologically irreducible representation of a central extension of $\ggG^X$ with central term $c(a,b)$ given by Eq. (\ref{c}).
\end{theoremdemo}
\begin{pf}
The proofs that $L^Y(a)$ satisfy Eq. (\ref{extended}) and are antihermitian require a straightforward calculation. So does the fact that the bilinear form $c(a,b)$ (Eq. (\ref{c})) is indeed a 2-cocycle - the antisymmetry is explicit, while the cocycle condition 
\[
c([a,b],c) + c([c,a],b) + c([b,c],a)=0
\]
follows from the invariance of the inner product $\langle\ ,\ \rangle$ on $\ggG$.
In order to prove irreducibility, we notice that the operators $L^Z(a)$ in Eq. (\ref{realrep}) are real, when $Z$ is a real field, i.e preserve separately the real and imaginary parts of our representation space. Therefore the operators $L^Z(a)$ define an irreducible  representation on the real part of our complex representation space. For a complex vector field $Y=Z_1+iZ_2$ we have $L^Y(a)=L^{Z_1}(a)+iL^{Z_2}(a)$. If we consider our complex representation space over the reals, then it is (real-)isomorphic to a direct sum of two copies of its real  part, i.e.
$V\cong V_{\bBr}\oplus V_{\bBr}$ and the operators $L^Y(a)$ are written in a block form:
\[
L^Y(a)=\left(\begin{array}{cc}
L^{Z_1}(a) & -L^{Z_2}(a)  \\
L^{Z_2}(a) & L^{Z_1}(a)\end{array}\right)\ :V_{\bBr}\oplus V_{\bBr} \rightarrow V_{\bBr}\oplus V_{\bBr}\ \ .
\]
Suppose that there is a nontrivial invariant subspace $U\subset V$. Then $U\cong U_{\bBr}\oplus U_{\bBr}$ and any element of 
$U_{\bBr}\oplus \{0\}$ should remain in $U_{\bBr}\oplus U_{\bBr}$ under the action of any $L^Y(a)$, i.e. 
\[
L^Y(a)(U_{\bBr}\oplus \{0\})\subset U_{\bBr}\oplus U_{\bBr}\ \ ,\ \ \ \forall a\in \ggG^X
\]
But that would imply that $L^{Z_1}(a)(U_{\bBr})\subset U_{\bBr}\ ,\forall a$, which contradicts the irreducibility of the representation on $V_{\bBr}$. This ends the proof.
\end{pf}
It is a simple matter to extend the above results to representations of central extensions of $\ggG_{\bBc}^X$. The operators must be defined as follows:
\begin{equation}
L^Y(a):=L_0(a)+B^*(Ya)-B(Y{\overline{a}})\,, \qquad a\in\ggG_{\bBc}^X\,. \label{complex}
\end{equation}
\newtheorem{theoremdemo2}{CORROLARY}
\begin{theoremdemo2}
The operators defined by {\rm Eq. (\ref{complex})} with $Y=Z_1+iZ_2$
give a topologically irreducible representation of a central extension of 
$\ggG_{\bBc}^X$ with central term given by: 
\begin{equation}
c(a,b)=-2i\left[(Z_1\overline{b},Z_2a)-(Z_1\overline{a},Z_2b)\right] \ .\label{ccomp}
\end{equation}
\end{theoremdemo2}
\end{section}
\begin{section}{Etingof--Frenkel current algebras and their representations}
Let $X$ be a compact Riemann surface of genus $g$. The Etingof--Frenkel current algebras are central extensions of the algebras $\ggG_{\bBc}^X$, defined as follows \cite{Etingof}:
Denote by $H_X$ the $g$-dimensional complex space of holomorphic differentials
on $X$ and by $H^*_X$ its dual. For every $a,b\in\ggG^X_{\bBc}$ and
$\alpha\in H_X$ one gets a number:
\begin{equation}
\Omega(a,b)(\alpha):=\int_X\langle a\,,\,db\rangle\wedge\alpha\,. \label{omega}
\end{equation}
Thus $\Omega(a,b)\in H^*_X$.
Using the fact that $\alpha$ is holomorphic and hence closed,
one deduces that $\Omega(a,b)(\alpha)$ is antisymmetric in $a$ and $b$ and
due to the invariance of the inner product $\langle\ ,\ \rangle$
it satisfies the cocycle condition.  Thus $\Omega$ defines a $g$-dimensional
central extension of $\ggG^X_{\bBc}$. This extension
has the important property that it is integrable to a
central extension by a $g$-dimensional complex torus of the corresponding group
$G^X_{\bBc}$ \cite{Etingof}.\par
   Rewriting the central term in Eq. (\ref{ccomp}), we can make it look similar to the one in Eq. (\ref{omega}). This is most easily done in local coordinates. Without loss of generality we may assume that the Riemann-Lebesgue measure in Eq. (\ref{scalar}) is just $dxdy$, since any factor, multiplying it can be absorbed into the vector field $Z$, leading, according to Prop. 1, to an equivalent representation. In local coordinates the two vector fields $Z_1$ and $Z_2$ can be written as
$Z_i=f_i{\partial\over\partial x}+g_i{\partial\over\partial y}\ ,\ i=1,2$. Performing integration by parts, keeping in mind that the vector fields are allowed to have jumps, we get:
\begin{equation}
c(a,b)=-2i\int_X\langle a\,,\,{\partial b\over\partial x}{\partial\varphi\over\partial y}-{\partial b\over\partial y}{\partial\varphi\over\partial x}\rangle\  dxdy\ +\ \hbox{boundary
terms}\,,
\end{equation}
where $\varphi:=f_1g_2-f_2g_1$.
The same manipulation, explained in invariant notations is as follows:
Let $\omega_{1,2}$ be the (piecewise smooth)
one-forms, corresponding to the vector fields $Z_{1,2}$ in the identification
of the tangent and cotangent spaces via the Riemannian metric.
Let $\varphi:=(\omega_1\wedge\omega_2)^*$ be the 0-form, which is Hodge dual
to $\omega_1\wedge\omega_2$. Then
\begin{equation}
c(a,b)=-2i\int_X\langle a\,,\,db\rangle\wedge d\varphi\ +\ \hbox{boundary
terms}.\label{myomega}
\end{equation}
Comparing Eqs. (\ref{omega}) and (\ref{myomega}), we notice several differences. The central term of Etingof and Frenkel defines a complex $g$-dimensional extension, while the ones in Eq. (\ref{myomega}), having fixed the vector field $Z$, are (real) one-dimensional. However, in any irreducible representation, the center should be projected to
$\bBC$, which should then be represented as a
multiplication operator. This is achieved by fixing a holomorphic one-form
$\alpha$ in (\ref{omega}). This will be (part of) the data, determining the
representation. Second, the function $d\varphi$ in (\ref{myomega}) is real, while $\alpha$ in (\ref{omega}) is complex. This can be taken care of by the following generalization:
We choose two different complex vector fields
$Y^L$ and $Y^R$, which have the same real part, (the latter requirement is needed
in order to ensure irreducibility.) We define new operators $L(a)$ (suppressing the superscript in $L(a)$, indicating its dependence on $Y^L$ and $Y^R$):
\begin{equation}
L(a)=L_0(a)+B^*(Y^Ra)-B(Y^L\overline{a})\,. \label{twovect}
\end{equation}
These give a representation of a central extension of $\ggG^X_{\bBc}$ with the
following central term:
\begin{equation}
c(a,b)=\left[(Y^L\overline{b},Y^Ra)-(Y^L\overline{a},Y^Rb)\right]\,.\label{ctwo}
\end{equation}
Rewriting the central term (\ref{ctwo}) we now get the same expression as in (\ref{myomega}) with a complex function
$\varphi:=f^R\overline{g^L}-\overline{f^L}g^R$ with $f^{L,R}$ and
$g^{L,R}$ being the corresponding complex coefficient functions of the vector
fields $Y^{L,R}$. Since $\varphi$ must be holomorphic, we see that the coefficient functions of $Y^R$ must be holomorphic, while those of $Y^L$ --- antiholomorphic. We should point out that the fields $Y^{L,R}$ cannot be chosen antiholomorphic and holomorphic, respectively. This would lead to $\varphi=0$.\par
   An essential difference between the central terms in (\ref{omega}) and (\ref{myomega}) is that the one-form $d\varphi$ in 
(\ref{myomega}) is (explicitly) exact, while the one-form $\alpha$ in (\ref{omega}) is never exact, since there are no exact holomorphic differentials on a compact Riemann surface (see, e.g. \cite{Kra}). At this point we recall that
the vector fields $Y^{L,R}$, and hence $\varphi$, can have jumps. Let us therefore choose $2g$ cycles $\gamma_j,
 \ j=1,2,...,2g$ on $X$, forming a canonical homology basis of $X$ (i.e. each $\gamma_j, j=1,...,g$ intersects once  $\gamma_{j+g}$ and does not intersect any other cycle). On the simply connected domain $X-\cup_j\gamma_j$ we fix a point $P_0$ and define
\begin{equation}
\varphi(P)=\int_{P_0}^P \alpha\ \label{alpha}
\end{equation}
For a fixed holomorphic differential $\alpha$, let us  denote by $c_j$ its $j$-th period, i.e. $c_j=\int_{\gamma_j}\alpha$ and set
$\hat{c}_j:=c_{j+g},\  j=1,...,g,\ \ \hat{c}_j:=c_{j-g},\  j=g+1,...,2g$.
 Then $\varphi$ has a jump equal to $-\hat{c}_j$ everywhere along the cycle $\gamma_{j}$ (the sign of the jump along a cycle is with respect to the positive direction of its dual cycle).
Now the operators $L(a)$ defined by (\ref{twovect}) with $\varphi$ given by (\ref{alpha}) will give a representation with
a central term
\begin{equation}
c(a,b)=\int_X\langle a\,,\,db\rangle\wedge \alpha\ +\
\sum_{j=1}^{2g}\hat{c}_j\int_{\gamma_j}\langle a\,,\,db\rangle\,.\label{bterm}
\end{equation}
Let $J=\{k\in\{1,..., 2g\}\,; \hat{c}_k\not= 0\}$. For each $j\in J$ choose an irreducible representation $l_j$ of 
the (unique) central extension of the loop algebra $\ggG^{S^1}_{\bBc}$
 with central charge equal to $-\hat{c}_j$. \par
Recall that the loop algebra $\ggG^{S^1}_{\bBc}$ has, up to a constant, a unique central extension (see, e.g. \cite{Pressley}), determined by the 2-cocycle
\begin{equation}
c(a,b)=\int_{S^1}\langle a\,,\,db\rangle\,\ \ a,b\in \ggG^{S^1}_{\bBc}\,.  \label{circle}
\end{equation}
Every centrally extended loop algebra $\ggG^{S^1}_{\bBc}$ contains as a dense subalgebra the algebra of Fourier polynomials. The latter  is, according to the classification of Kac \cite{Kac}, an untwisted affine algebra (after adjoining a one-dimensional space of derivations). If $\ggB$ is such an algebra of Fourier polynomials, then every irreducible highest--weight (or lowest--weight) module over $\ggB$ gives rise to an irreducible representation of the centrally extended loop algebra, with the same central charge.\par
For each cycle $\gamma_j$ there is an obvious irreducible representation of the extension of $\ggG^X_{\bBc}$ with central term which looks as in (\ref{circle}) obtained by restricting any $a\in \ggG^X_{\bBc}$ to the cycle $\gamma_j$ and then representing it via $l_j$. Thus, by a slight abuse of notation, we can define the following representations:
\begin{equation}
{\cal L}:=L\otimes\left(\mathop{\otimes} \limits_{j\in J}l_j\right)\ .  \label{last}
\end{equation}
Here the tensor products are in the sense of tensoring of representations of
a Lie algebra. Strictly speaking we are  tensoring representations of different algebras
but as the difference is just in the central terms, no difficulty arises. \par
We can now prove the main result of our paper:
\begin{theoremdemo}\label{Prop3}
Fixing a holomorphic differential $\alpha$ on the Riemann surface $X$, the representations ${\cal L}$
defined by {\rm Eq. (\ref{last})} are topologically irreducible representations of the Etingof-Frenkel algebra with central extension given by {\rm Eq. (\ref{omega})}.
\end{theoremdemo}
\begin{pf}
The fact that these are representations of the Etingof-Frenkel algebras is obvious from the construction. The added factors $l_j$ contribute to the central charge exactly the necessary terms to cancel the boundary terms in Eq. (\ref{bterm}) and what remains is the central term of Eq. (\ref{omega}). It is less trivial to show irreducibility. If we denote by $X_0$ the open submanifold $X-\cup_j\gamma_j$, the main idea is to show that we can essentially treat the algebra $\ggG^{X_0}_{\bBc}$ independently from the algebra of functions on $\cup_j\gamma_j$ with values in $\ggG_{\bBc}$. In addition the module $\mathop{\otimes}_{j\in J}l_j$ could in principle be reducible under the action of the latter algebra, since
the cycles $\gamma_j$ are not disjoint (they are linked two by two) and therefore our algebra is a subalgebra of the tensor product of loop algebras. Topological irreducibility will follow easily from the following two lemmas.
\end{pf}
\newtheorem{theoremdemo1}{Lemma}
\begin{theoremdemo1}
Let $X_0\subset X$ be an open submanifold of full measure. The representation defined by {\rm Eq. (\ref{realrep})} is topologically irreducible under the action of the subalgebra $C^{\infty}_0(X_0, \ggG)\subset \ggG^X$ consisting of smooth functions with compact support inside $X_0$. Similarly, the representation of the central extension of $\ggG^X_{\bBc}$ given by {\rm Eq. (\ref{complex})} remains topologically irreducible under the action of the subalgebra 
$C^{\infty}_0(X_0, \ggG_{\bBc})$
\end{theoremdemo1}
\begin{pf}
The proof of this fact (at the group level) is contained in \cite{Gelfand} and is also discussed in \cite{Donkov}. At the Lie algebra level it can be easily deduced from the explicit action of the operators $L_0(a)$, $B^*(Za)$ and $B(Z{a})$ and the fact that any function in the one--particle
Hilbert space $L^2(X,dx;\ggG_{\bBc})$ can be approximated (in the $L^2$ sense) by smooth functions with compact support inside $X_0$.
\end{pf}
\begin{theoremdemo1}
Let $\ggB^0_1$ and $\ggB^0_2$ be two copies of the algebra of Fourier polynomials on the unit circle $S^1$ with values in $\ggG_{\bBc}$, i.e. 
\[
\ggB^0_j\cong \bBC[e^{ 2\pi i t_j}, e^{- 2\pi i t_j}]   \otimes\ggG_{\bBc}\,,\  
\ t_j \hbox{ \rm -- formal variables}\,,\ \ j=1,2.
\]
Let $\ggB_{1,2}$ be their central extensions, determined by the 2-cocycle in {\rm Eq. (\ref{circle})}. Denote by $\ggD$ the
subalgebra of $\ggB_1\otimes\ggB_2$ generated by pairs of Fourier polynomials, coinciding when evaluated at 
$t_1=t_2=0$. 
Take two irreducible highest--weight modules $l_{1,2}$ over $\ggB_{1,2}$, respectively. Then the tensor product 
$l_1\otimes l_2$ (which is obviously an irreducible module over $\ggB_1\otimes\ggB_2$) remains irreducible under the action of  $\ggD$. (The statement remains true for two lowest--weight modules or a product of a highest--weight and lowest--weight modules.)
\end{theoremdemo1}
\begin{pf}
It is enough to show that for any two vectors $u,v \in l_1\otimes l_2$ we can find an element of the universal enveloping algebra $U(\ggD)$ which sends $u$ to $v$. Let us first consider a homogeneous element $u=u_1\otimes u_2\in l_1\otimes l_2$. Since $l_1$ is an irreducible highest--weight module (with some highest weight $\Lambda_1$), there is an element 
$a\in U(\ggB_1)$, which sends $u_1$ to the highest--weight vector $u^\Lambda_1$. Note that we have the following identifications of vector spaces: 
\begin{eqnarray}
\nonumber U(\ggB_1)&=&U\left(\bBC[e^{ 2\pi it_1}, e^{- 2\pi it_1}]   \otimes \ggG_{\bBc}\oplus \bBC\right)\\
\nonumber &=&U\left(\bBC[e^{ 2\pi i t_1}, e^{- 2\pi i t_1}]\otimes \ggG_{\bBc}\right)   \\
\nonumber &=&\bBC[e^{ 2\pi i t_1}, e^{- 2\pi i t_1}]\otimes U\left(\ggG_{\bBc}\right)
\end{eqnarray}
The vector $u_2$ can be decomposed into a sum of weight vectors:
\[
u_2=\sum\limits_{\lambda\leq \Lambda_1}u_2^\lambda\ ,
\]
where only a finite number of elements of this sum are nonzero. Thus we can find big enough positive integer $n$, so that
\[
(e^{ 2\pi int_2}\otimes g)\, u_2=0\,,\ \  \forall g\in U(\ggG_{\bBc})\,.
\]
Let $g_a\in U(\ggG_{\bBc})$ be the image of $a\in U(\ggB_1)$ under the evaluation homomorphism (evaluating the Fourier polynomial at $t_1=0$). We choose an element $\hat{a}=e^{ 2\pi int_2}\otimes g_a\in U(\ggB_2)$. Then
\[
(a\otimes\hat{a})(u_1\otimes u_2)=(a\, u_1)\otimes u_2+a\otimes(\hat{a}\, u_2)=u^{\Lambda_1}\otimes u_2
+0\,.
\]
Note that the element $a\otimes\hat{a}$ actually belongs to the universal enveloping algebra of the subalgebra $\ggD$.
  Similarly, we can find suitable element of $U(\ggD)$ which will send $u^{\Lambda_1}\otimes u_2$ to $u^{\Lambda_1}\otimes u^{\Lambda_2}$.\par
Let now $\sum_j u_1^j\otimes u_2^j$ be a general element of $l_1\otimes l_2$. Again, choosing appropriately
$(a\otimes\hat{a})\in U(\ggD)$ as before with big enough positive $n$, we can bring some of the elements of the sum to the form $\beta_j u^{\Lambda_1}\otimes u_2^j$, while others will turn to zero. Grouping the remaining nonzero elements we get $u^{\Lambda_1}\otimes u_2$ for some $u_2$. On the next step we can get  
$u^{\Lambda_1}\otimes u^{\Lambda_2}$. \par
  Then, by choosing appropriately $(b\otimes\hat{b})\in U(\ggD)$ the element $u^{\Lambda_1}\otimes u^{\Lambda_2}$ is first sent to $u^{\Lambda_1}\otimes v_2$ and then to $v_1\otimes v_2$. Since any $v\in l_1\otimes l_2$ is a finite sum of such homogeneous vectors, we have shown that an appropriate linear combination of elements of $U(\ggD)$ will send $u$ to $v$. This completes the proof.
\end{pf}
\end{section}
\begin{section}{Examples. The Problem of Integrability}
There seems to be an intriguing connection between the theory of Riemann surfaces and the question of possibility or impossibility to obtain integrable irreducible representations using the scheme described.
Unfortunately, it appears impossible to get integrable representations using just highest or lowest weight (integrable) modules $l_j$ in Eq. (\ref{last}). \par
 Let us fix a canonical homology basis $\{\gamma_j\}_{j=1}^{2g}$ on $X$, and find the so-called dual basis $\{\alpha_i\}_{i=1}^n$ of holomorphic differentials, defined by the requirement that their first $g$ periods (commonly named $a$-periods) are
\[\int\limits_{\gamma_j}\alpha_i=\delta_{ij}\,,\ \ i,j=1,2,...,g\ .\]
The remaining $g$ periods of $\alpha_i$, the so-called $b$-periods,  form a symmetric $g\times g$ matrix with positive-definite immaginary part. This is the period matrix of $X$ - a very important invariant of the complex structure.\par
A  sufficient condition for integrability of the representation  (\ref{last}) would be certain integrality condition for the holomorphic differential $\alpha$ (integral linear combination of the basis elements $\alpha_i$) plus integrability of the modules $l_j$ in  
(\ref{last}). Among all modules $l_j$ coming from highest-- (lowest--)weight modules over Kac--Moody algebras, integrable are only the ones, for which the central term is a positive (negative) immaginary integer \cite{Kac}.
Because the $a$-periods and $b$-periods appear in our setting as the central charges of the modules $l_j$ that we have to choose in Eq. (\ref{last}), we might be tempted to try to choose a complex structure on $X$ and a holomorphic differential $\alpha$ in such a way, that all periods are immaginary integers. This is impossible due to the positive definiteness of the immaginary part of the period matrix. The best we can hope for is to have purely real integer $a$-periods and purely immaginary integer $b$-periods. This will become clear from the two examples we consider.\par
Consider first the torus: $X=\bBC/\Lambda_{1,i}$ (the complex plain factorized by the lattice, generated by $1$ and $i$). The image of the real axis is an $a$-cycle, and of the imaginary axis --- a $b$-cycle. The holomorphic differential $dz$ forms the dual basis. We have $c_1=\int_adz=1$ and $c_2=\int_bdz =i$, ($\hat{c}_1=i$, $\hat{c}_2=1$). Thus we need to add in (\ref{last}) as factors an (integrable) lowest--weight module with central charge $-i$ and a (nonintegrable) lowest--weight module with central charge $-1$.\par
\begin{figure} [h!]
\centerline{
\includegraphics[scale=0.65]{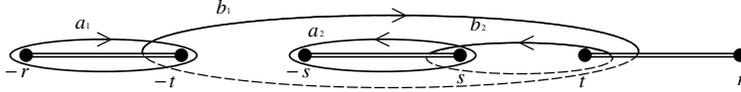}  }                              
\caption{A hyperelliptic surface of genus 2 } 
\label{hyper}
\end{figure}
A less trivial example is a hyperelliptic surface of genus 2, (Figure \ref{hyper}). This is a two-sheeted cover of the Riemann sphere 
$\bBC\cup \{\infty \}$ with six branch points $0<s< t< r$ and
$-r,-t,-s$ on the real line. Make three cuts, as shown on the picture, glueing the upper banks of the front plane with the lower banks of the back plane and vice-versa. This is topologically a sphere with two handles. The two $a$-cycles and the two $b$-cycles can be chosen as depicted and constitute a canonical basis (note that the dashed part of the $b$-cycles indicates that they pass through the cuts onto the lower plane. The following is a holomorphic differential on $X$ \cite{Kra}
\[
\alpha:={z\,dz\over\sqrt{(z^2-s^2)(z^2-t^2)(z^2-r^2)}}
\]
Let's calculate the periods of the differential $\alpha$. Denote by $R(x;s,t,r)$ the function 
$x[(x^2-s^2)(x^2-t^2)(x^2-r^2)]^{-{1\over 2}}$. Note that $R(x;s,t,r)$ is odd.
Keeping in mind that for $z$ immediately above or below a cut, $\sqrt{z}$ is purely imaginary, while for $z$ on the real line, where there are no cuts, $\sqrt{z}$ is real, we obtain:
\[
\int\limits_{a_1}\alpha=-2i\int\limits_{-r}^{-t}R(x;s,t,r)dx=2i\int\limits_{t}^{r}R(x;s,t,r)dx
\]
\[
\int\limits_{b_2}\alpha=-2\int\limits_{s}^{t}R(x;s,t,r)dx
\]
The other two periods are zero since the intervals of integration are symmetric w.r.to zero.
The two nonzero integrals above can be expressed in terms of elliptic integrals of the first kind.
Since we have plenty of freedom to vary $s$, $r$ and $t$, we can make the two periods mutually rational. Note that the first is purely immaginary, the second is real. Then an appropriate multiple of $\alpha$ will have the property that its $a$-periods are immaginary integers and the $b$-periods are real integers.\par
It becomes obvious that one has to go beyond highest and lowest weight modules if one hopes to construct integrable irreducible representations.\par
\end{section}\bigskip
{\bf Acknowledgements}\par\smallskip\noindent
The author would like to thank Emil Horozov and Vassil Tsanov for helpful discussions.

\end{article}
\end{document}